\documentclass{amsart}
\usepackage{amscd}
\usepackage{amsmath}
\usepackage{amsxtra}
\usepackage{amsfonts}
\usepackage{amssymb}

\newtheorem{theorem}{Theorem}[section]
\newtheorem{corollary}[theorem]{Corollary}
\newtheorem{lemma}[theorem]{Lemma}
\theoremstyle{definition}
\newtheorem{remark}[theorem]{Remark}
\newtheorem{example}[theorem]{Example}
\theoremstyle{remark}
\newtheorem*{acknowledgements}{Acknowledgements}

\numberwithin{equation}{section}

\newcommand{\ip}[2]{\left\langle\,#1\mid#2\,\right\rangle}
\newcommand{\cj}[1]{\overline{#1}}
\newcommand{\bz}{\mathbb{Z}}
\newcommand{\br}{\mathbb{R}}
\newcommand{\bc}{\mathbb{C}}
\newcommand{\bt}{\mathbb{T}}
\newcommand{\bn}{\mathbb{N}}
\newbox\frogdown
\long\def\hookdownarrow{\setbox\frogdown=\hbox to 0pt{\hss $\displaystyle \downarrow $\hss }\vrule height0pt width0pt depth1.5\ht\frogdown
\setlength{\unitlength}{0.4pt}
\begin{picture}(10,5)
\put(5,0){\oval(10,10)[t]}
\end{picture}
\lower\ht\frogdown\box\frogdown}

\begin{document}
\title[Wavelet constructions in non-linear dynamics]{Wavelet constructions in non-linear dynamics}
\author{Dorin Ervin Dutkay and Palle E.T. Jorgensen}
\address{Department of Mathematics\\
The University of Iowa\\
14 MacLean Hall\\
Iowa City, IA 52242-1419\\
U.S.A.}
\email{Dorin Ervin Dutkay: ddutkay@math.rutgers.edu}\
\email{Palle E.T. Jorgensen: jorgen@math.uiowa.edu}
\begin{abstract}
We construct certain Hilbert spaces associated with a class of
non-linear dynamical systems $X$. These are systems which arise from a
generalized self-similarity, and an iterated substitution. We show that when
a weight function $W$ on $X$ is given, then we may construct associated Hilbert
spaces $H(W)$ of $L^2$-martingales which have wavelet bases.
\end{abstract}
\maketitle

\section{Introduction}
    A particularly productive approach to the construction of wavelet bases
in $L^2(\br)$ is based on the notion from optics of \emph{resolution}, which translates
into scales of nested Hilbert subspaces $V_n$, $n \in \bz$, in $L^2(\br)$ such that the
intersection is $\{0\}$ and the union is dense. Moreover the operation of dyadic
scaling transforms each $V_n$ to the next $V_{n+1}$. This is called a
multiresolution approach to wavelets, see \cite{Dau92}, and it is based on the
interplay between the two abelian groups $\bt$ (the circle group${}={}$one-torus), and
$\br$ (the real line), with $\bt$ representing a period interval placed on the line
$\br$. This paper is based on the observation that multiresolutions really are
martingales; and by exploiting this fact, we are able to adapt the geometric
idea of subspaces ($V_n$) to nonlinear dynamics in a variety of applications
where such a pair of groups is not available, but instead there is a single
endomorphism on a compact space $X$ which defines a certain self-similarity
mirroring the more familiar scale-similarity that is so powerful in wavelet
analysis.
\par
The words `non-linear' and `wavelets' in the title beg two questions: (1)
``What is the Hilbert space?''  (2) If `non-linear', then there must be a
substitute for the duality between the operators of translation and
multiplication !?
\par
We will address the questions in the announcement below, while giving
answers with full proofs in forthcoming papers \cite{DJ03,DJ04a,DJ04b,DJ04c}
\par
We have in mind three classes of examples: (1) The state space $X$ for a
sub-shift system in symbolic dynamics; (2) affine iterated function systems
based on a fixed expansive scaling matrix; and (3) the complex iteration
systems which generate Julia sets $X$ in the Riemann sphere. If $r(z)$ is a
 rational function, set $r^n=\underbrace{r\circ\dots\circ r}_{n\text{ times}}$.
 Then the Julia set $X=X(r)$ is the complement of the largest open subset of
 $\bc$ where $r^n$ is a normal family.
\par
The first step to construct a Hilbert space in any of these three classes
of examples is the identification of the appropriate covariant measures on
$X$, and the second step is the construction of a certain lifting from $X$ to a
suitable space $X_\infty$ of discrete paths which originate in $X$.
Specifically,
$X_\infty:=\{\,(x_n)\in X^{\bn}\mid r(x_{n+1})=x_n,\;n\in\bn\,\}$.
      We begin (in Section \ref{LNLMart}) by recalling how the familiar wavelet
construction may be obtained from $X_\infty$ in the special case when $X = \bt
= \{\,z \in \bc\mid|z| = 1\,\}$, and $r(z) = z^2$. While the real line $\br$ is not
$\bt_\infty$, we may still build the Hilbert space $L^2(\br)$ with its
multiresolution wavelets from an inductive and isometric procedure based on
$\bt_\infty$. The intuitive idea is to get $\br$ in the limit by successive
doubling of periods. But we use Hilbert spaces, and the idea is outlined in
(\ref{eqint3}). It is further based on martingales. To highlight the distinction
between the compact group $G = \bt_\infty$ (called a solenoid) and the reals
$\br$, contrast the following two familiar short exact sequences of abelian
groups:
$0\rightarrow C\rightarrow G\rightarrow\bt\rightarrow0$,
where $C$ is the Cantor group, on the one hand; and
$0\rightarrow\bz\hookrightarrow\br\rightarrow\bt\rightarrow0$, on the other.
So while the groups $G$ and $\br$ are
very different, the use of $G$ helps us build the Hilbert space $L^2(\br)$, but
within the category of isometries in Hilbert space. Moreover, the $X_\infty$
viewpoint is useful in important applications outside the group context.

We now introduce weight functions which are determined from our
candidate for low-pass filters $m$. In our analysis, a low-pass filter will be
a function $m$ on $X$ which satisfies axioms that generalize those which are
known for standard wavelets in $L^2(\br)$. Similarly the relation between the
circle $\bt$, and the real line $\br$, for familiar wavelets is mirrored in our
lifting from $X$ to the space $X_\infty$.
\par
A second feature in our analysis is a certain Perron-Frobenius-Ruelle
operator $R=R_W$ associated with a non-negative weight function $W$ which may
be taken to be the absolute value-square of $m$.
\par
In addressing the second question, we give up the analogue of
translations, but for the case of Julia sets $X$, we work instead
with the operator of multiplication by the variable $z$ restricted
to $X$. Our Hilbert space will be a Hilbert space built on
martingales on $X_\infty$, and the unitary operator $U$ which
corresponds to the familiar dyadic, or $N$-adic scaling for
familiar wavelets will simply be the substitution of our system.
This operator $U$ will scale between levels of our discrete
$L^2$-martingales. Hence we arrive at a class of generalized
multiresolutions, or multi-wavelets. In this setting, we are able
to prove a version of the dimension consistency formula of
Baggett-Merrill et al.\ \cite{BaMe99}, and to construct our wavelets.
\par
A new feature of our analysis is a dichotomy between properties of our
low-pass filters: The low/high-pass conditions in our context will refer to
either a finite set of points or to a singular measure with full support. In
the case of the middle-third Cantor set, this measure will be a classical
infinite product Riesz measure.
\par
What the non-linear settings have in common with the classical linear
cases (e.g., wavelets) is a certain dichotomy. The non-linear systems we
consider carry a certain strongly invariant measure. The support of this
invariant measure must either be full, or it must be finite. In the case
when the support is full we can expect the same rigidity as is known in
familiar wavelet constructions, namely that the filters giving wavelets in
the same space have equal absolute value. When the support is finite, we
construct some scaling function as an infinite product, and we show that our
generalized ``low-pass'' condition implies the existence of a wavelet basis
construction is in a fixed Hilbert space, avoiding redundant multiplicity.
We offer results on the nature of the Hilbert space which carries wavelet
bases for both cases in our dichotomy: full support, and finite support. In
either case, our construction uses a family of discrete time $L^2$-martingales.
 We establish a geometric setting which does not depend on the
generalized wavelet filter. We then give conditions on two filters which
imply isomorphism of the corresponding wavelet systems when the support of
our invariant measure is some given finite cycle.
\par
In the case of wavelets, \cite{BDP04}, it is $L^2(\br)$, or a finite
direct sum of $L^2(\br)$ with itself. In that case the complex function is
$z\rightarrow z^N$, and the Julia set is the circle $\bt$.
\par
Due to results of Bramson and Kalikow \cite{BK93} the strongly invariant measure
 may be non-unique, even in the case of a full shift, if the weight $W$ is just
continuous.

\section{\label{LNLMart}Linear vs non-linear: martingales}
In this paper we study the problem of inducing operators on
Hilbert space from non-invertible transformations on compact
metric spaces. The operators, or representations must satisfy
relations which mirror properties of the given point
transformations.
\par
While our setup allows a rather general formulation in the context
of $C^*$-algebras, we will emphasize the case of induction from an
abelian $C^*$-algebra. Hence, we will stress the special case when
$X$ is a given compact metric space, and $r\colon X\rightarrow X$ is a
finite-to-one mapping of $X$ onto $X$. Several of our results are
in the measurable category; and in particular we are not assuming
continuity of $r$, or any contractivity properties.
\par
If $r\colon X\rightarrow X$ (onto) is given, we let $X_\infty$ denote the projective limit
\begin{equation}
\begin{array}{cccccccc}
X&\overset{\textstyle r}{\longleftarrow}&
X&\overset{\textstyle r}{\longleftarrow}&
X&\cdots&\longleftarrow&X_\infty
\,,
\end{array}
\label{Xinfty}
\end{equation}
i.e.,
$
X_\infty:=\{\,(x_n)_{n\in\mathbb{N}_0}\mid r(x_{n+1})=x_n\,\}
$.
Then $r$ extends to an automorphism $\hat r$ determined by
\[
\begin{array}{ccc}
\phantom{X}\llap{$X_\infty$} & 
\overset{\textstyle \hat{r}}{\longrightarrow}
& \rlap{$X_\infty$}\phantom{X} \\ 
\makebox[0pt]{\hss \rule[-8pt]{0.4pt}{14pt}\hss }\setbox\frogdown=\hbox to 0pt{\hss $\displaystyle \downarrow $\hss }\lower\ht\frogdown\box\frogdown
&  & 
\makebox[0pt]{\hss \rule[-8pt]{0.4pt}{14pt}\hss }\setbox\frogdown=\hbox to 0pt{\hss $\displaystyle \downarrow $\hss }\lower\ht\frogdown\box\frogdown\rule[-12pt]{0pt}{20pt} \\ 
X & 
\overset{\smash{\textstyle r}}{\longrightarrow}
& X
\end{array}
\]

\subsection{Wavelets}
\par Our results will apply to wavelets. In the theory
of multiresolution wavelets, the problem is to construct a special
basis in the Hilbert space $L^2(\br^d)$ from a set of numbers
$a_n$, $n\in\bz^d$.
\par
The starting point is the scaling identity
\begin{equation}\label{eqintsc1}
\varphi(t)=N^{1/2}\sum_{n\in\bz^d}a_n\varphi(At-n),\quad(t\in\br^d),
\end{equation}
where $A$ is a $d$ by $d$ matrix over $\bz$, with eigenvalues
$|\lambda|> 1$, and $N=|\operatorname*{det} A|$, and where $\varphi$ is a
function in $L^2(\br^d)$.
\par
The first problem is to determine when (\ref{eqintsc1}) has a
solution in $L^2(\br^d)$, and to establish how these solutions
({\it scaling functions}) depend on the coefficients $a_n$.
\par
When the Fourier transform is applied, we get the equivalent
formulation,
\begin{equation}\label{eqintsc2}
\hat\varphi(x)=N^{-1/2}m_0({A^{tr}}^{-1}x)\hat\varphi({A^{tr}}^{-1}x),
\end{equation}
where $\hat\varphi$ denotes the Fourier transform,
\[\hat\varphi(x)=\int_{\br^d}e^{-i2\pi x\cdot t}\varphi(t)\,dt\]
and where now $m_0$ is a function on the torus
\[
\bt^d=\{\,z=(z_1,\dots,z_d)\in\bc^d\mid
 |z_j| =1, \;1\leq j\leq d\,\}=\br^d/\bz^d
,
\]
i.e.,
$m_0(z)=\sum_{n\in\bz^d}a_nz^n=\sum_{n\in\bz^d}a_ne^{-i2\pi
n\cdot x}$ ($z=e^{-2\pi i\cdot x}$). The duality between the
compact group $\bt^d$ and the lattice $\bz^d$ is given by
$\ip{z}{n}=z^n=z_1^{n_1}\cdots z_d^{n_d}$,
$z=(z_1,\dots,z_d)$,
$n=(n_1,\dots,n_d)$.
\par In this case, matrix multiplication $x\mapsto A x$ on $\br^d$ passes to the quotient $\br^d/\bz^d$, and we get an $N$-to-one
mapping $x\mapsto Ax\bmod\bz^d$, which we denote by $r_A$.
\par
The function $m_0$ is called a low pass filter, and it is chosen
such that the operator $S=S_{m_0}$ given by
$(Sf)(z)=m_0(z)f(Az)$
is an isometry on $H_0=L^2(\bt^d$, Haar measure$)$. Moreover,
$L^\infty(\bt^d)$ acts as multiplication operators on $H_0$. If
$g\in L^\infty(\bt)$,
$(M(g)f)(z)=g(z)f(z)$,
and
\begin{equation}\label{eqintcov1}
SM(g)=M(g(A\,\cdot\,))S
\end{equation}
A main problem is the extension of this covariance relation
(\ref{eqintcov1}) to a bigger Hilbert space $H_0\rightarrow
H_{\operatorname*{ext}}$, $S\rightarrow S_{\operatorname*{ext}}$,
such that $S_{\operatorname*{ext}}$ is unitary
in $H_{\operatorname*{ext}}$. We now sketch briefly this extension in some
concrete cases of interest.
\par
We construct a sequence of measures
$\omega_0,\omega_1,\dots$ on $\bt^d$ such that
$L^2(\bt^d,\omega_0)\simeq H_0$, and such that there are natural
isometric embeddings
\begin{equation}\label{eqint3}L^2(\bt^d,\omega_n)\hookrightarrow
L^2(\bt^d,\omega_{n+1}),\quad f\mapsto f\circ r_A.\end{equation}
The limit in (\ref{eqint3}) defines a {\it martingale Hilbert
space} $\mathcal{H}$ in such a way that the norm of the
$L^2$-martingale $f$ is
$\|f\|^2=\lim_{n\rightarrow\infty}\|P_nf\|^2_{L^2(\bt^d,\omega_n)}$.
In \cite{DJ04a}, we also prove an associated  pointwise a.e.
convergence theorem.
\par
If $\Psi\colon  L^2(\bt^d,\omega_n)\rightarrow L^2(\br^d)$ is defined by
$\Psi\colon  f_n\mapsto f_n(A^{-n}x)\hat\varphi(x)$,
then $\Psi$ is an isometry of $L^2(\bt^d,\omega_n)$ into
$L^2(\br^d)$.
\par
Specifically,
\begin{equation}\label{eqint4}
\int_{\bt^d}|f_n|^2\,d\omega_n=\int_{\br^d}|f_n(A^{-n}x)\hat\varphi(x)|^2\,dx.
\end{equation}
As a result we have induced a system
$(r_A,\bt^d)\rightarrow (S_{m_0},L^2(\bt^d))\rightarrow
(U_A,L^2(\br^d))$ where
\begin{equation}\label{eqintdil}
(U_Af)(x)=N^{1/2}f(Ax)\quad(f\in L^2(\br^d),\;U_A
\text{ unitary}); \end{equation}
  the system is determined by the given filter function
$m_0$. It can be checked that $\Psi$ is
an isometry, and that
$U_AM(g)=M(g(A\,\cdot\,))U_A$
holds on $L^2(\br^d)$. Moreover $\Psi$ maps onto $L^2(\br^d)$ if
the function $m_0$ doesn't vanish on a subset of positive measure.
\par
In the case of wavelets, we ask for a wavelet basis in
$L^2(\br^d)$ which is consistent with a suitable resolution
subspace in $L^2(\br^d)$. Whether the basis is orthonormal, or
just a Parseval frame, it may be constructed from a system of
subband filters $m_i$, say with $N$ frequency bands. These filters
$m_i$ may be realized as functions on $X=\bt^d=\br^d/\bz^d$, the
$d$-torus. Typically the scaling operation is specified by a given
expansive integral $d$ by $d$ matrix $A$.
\par
 Let $N:= |\operatorname*{det} A|$. Pass $A$ to the quotient $X = \br^d/\bz^d$, and we get a mapping
$r$ of $X$ onto $X$ such that $\# r^{-1}(x) = N$ for all $x$ in
$X$, and the $N$ branches of the inverse are strictly contractive
in $X=\br^d/\bz^d$ if the eigenvalues of $A$ satisfy
$|\lambda|>1$.
\par
The subband filters $m_i$ are defined in terms of this map, $r_A$,
and the problem is now to realize the wavelet data in the Hilbert
space $L^2(\br^d)$ in such a way that $r=r_A \colon X\rightarrow X$
induces the unitary scaling operator $f\mapsto N^{1/2} f(A x)$ in
$L^2(\br^d)$, see (\ref{eqintdil}).
\subsection{Examples (Julia sets, subshifts)}
\par
In this paper we will show that this extension from spaces $X$,
with a finite-to-one mapping $r\colon X\rightarrow X$, to operator
systems may be done quite generally, to apply to the case when $X$
is a Julia set for a fixed rational function of a complex
variable, i.e., $r(z)=p_1(z)/p_2(z)$, with $p_1,p_2$ polynomials,
$z\in\bc$ and $N=\max(\operatorname*{deg}p_1,\operatorname*{deg}p_2)$. Then
$r\colon X(r)\rightarrow X(r))$ is $N$-to-$1$ except at the singular
points of $r$. Here $X(r)$ denotes the Julia set of $r$. \par
 It
also applies to shift-invariant spaces $X(A)$ when $A$ is a $0$--$1$
matrix and
$
X(A)=\left\{\,(x_i)\in\prod_{\mathbb{N}}\{1,\dots,N\}
\mid A(x_i,x_{i+1})=1\,\right\}
$,
while
$
r_A(x_1,x_2,\dots)=(x_2,x_3,\dots)
$
is the familiar subshift. Note that $r_A\colon X(A)\rightarrow X(A)$ is
onto iff every column in $A$ contains at least one entry $1$.
\subsection{Martingales}
\par
Part of the motivation for our paper derives from the papers by
Richard Gundy \cite{Gun00}, \cite{Gun04}, \cite{Gun99},
\cite{Gun66}. The second named author also acknowledges
enlightening discussions with R. Gundy. The fundamental idea in
these papers by Gundy et al.\ is that multiresolutions should be
understood as martingales in the sense of Doob
\cite{Doob1}, \cite{Doob2}, \cite{Doob3} and Neveu \cite{Neveu}. And
moreover that this is a natural viewpoint.
\par
  One substantial advantage of this viewpoint is that we are then able to
handle the construction of wavelets from subband filters that are
only assumed measurable, i.e., filters that fail to satisfy the
regularity conditions that are traditionally imposed in wavelet
analysis.
\par
  A second advantage is that the martingale approach applies to a number of
wavelet-like constructions completely outside the traditional
scope of wavelet analysis in the Hilbert space $L^2(\br^d)$. But
more importantly, the martingale tools apply even when the
operation of scaling doesn't take place in $\br^d$ at all, but
rather in a compact Julia set from complex dynamics; or the
scaling operation may be one of the shift in the subshift dynamics
that is understood from that thermodynamical formalism of David
Ruelle \cite{Rue89}.
\subsection{The general theory}
\par
 In each of the examples, we are faced with a given space $X$, and a
finite-to-one mapping $r\colon X\rightarrow X$. The space $X$ is
equipped with a suitable family of measures $\mu_h$, and the
$L^\infty$ functions on $X$ act by multiplication on the
corresponding $L^2$ spaces, $L^2(X,\mu_h)$. It is easy to see that
there are $L^2$ isometries which intertwine the multiplication
operators $M(g)$ and $M(g\circ r)$, as $g$ ranges over
$L^\infty(X)$. We have
\begin{equation}\label{eqintadd1}
\begin{array}{ccc}
\phantom{\mathcal{H}_{\operatorname*{ext}}}\llap{$L^2(X,\mu_h)$} & 
\overset{\textstyle S}{\longrightarrow}
& \rlap{$L^2(X,\mu_h)$}\phantom{\mathcal{H}_{\operatorname*{ext}}} \\ 
\hookdownarrow
&  & 
\hookdownarrow\rule[-10pt]{0pt}{20pt} \\ 
\mathcal{H}_{\operatorname*{ext}} & 
\overset{\smash{\textstyle U}}{\longrightarrow}
& \mathcal{H}_{\operatorname*{ext}}
\end{array}
\end{equation}
where the vertical maps are given by inclusions. Specifically,
\begin{equation}\label{eqintadd2}
SM(g)=M(g\circ r)S,\text{ and } UM(g)U^{-1}=M(g\circ r)
\end{equation}
\par
But for spectral-theoretic calculations, we need to have
representations of $M(g)$ and $M(g\circ r)$
that are
unitarily equivalent.
That is true in traditional wavelet applications, but the unitary
operator $U$ in (\ref{eqintadd2}) is not acting on $L^2(X,\mu_h)$.
Rather,
the unitary $U$ is acting by matrix scaling on a different
Hilbert space, namely $L^2(\br^d,\text{Lebesgue measure})$, with
the scaling operator $U_A$ from (\ref{eqintdil}).
\par In the other
applications, Julia set, and shift-spaces, we aim for a similar
construction. But in these other cases, it is not at all clear
what the Hilbert space corresponding to $L^2(\br^d)$, and the
corresponding unitary matrix scaling operator, should be.
\par
We provide two answers to this question, one at an abstract level,
and a second one which is a concrete function representation.
\par
At the abstract level, we show that the construction may be
accomplished in Hilbert spaces which serve as unitary dilations of
the initial structure, see (\ref{eqintadd1}). \par In the
concrete, we show that the extended Hilbert spaces may be taken as
Hilbert spaces of $L^2$-martingales on $X$. In fact, we present
these as Hilbert spaces 
of certain $L^2$-functions on the projective limit $X_\infty$ which we
outlined in (\ref{Xinfty}) above.
Only, now the relevant measures are part of this limit construction.
 This is analogous to the distinction between an
abstract spectral theorem on the one hand, and a concrete spectral
representation on the other. To know details about
multiplicities, and multiplicity functions, we need the latter.
\par
Our concrete version of the dilation Hilbert space
$\mathcal{H}_{\operatorname*{ext}}$ from (\ref{eqintadd1}) is then
\[\mathcal{H}_{\operatorname*{ext}}\simeq L^2(X_\infty,\hat\mu_h)\]
for a suitable measure $\hat\mu_h$ on $X_\infty$.
\section{\label{Fume}Functions and measures on $X$}
Consider
$X$ a compact metric space, 
with
$\mathfrak{B}$,
$r$,
$W$,
$\mu$ as follows:
$\mathfrak{B}=\mathfrak{B}(X)$ a Borel sigma-algebra of subsets of
$X$, 
$r\colon X\rightarrow X$ an onto, measurable map such that
$\# r^{-1}(x)<\infty$ for all $x\in X$, 
$W\colon X\rightarrow[0,\infty)$, 
$\mu$ a positive Borel measure
on $X$.
\subsection{Transformations
of functions and measures}
\begin{itemize}
\item Let $g\in L^\infty(X)$. Then
\begin{equation}\label{eqfm2_2_1}
M(g)f=gf
\end{equation}
is the multiplication operator on $L^\infty(X)$ or on
$L^2(X,\mu)$. \item Composition:
\begin{equation}\label{eqfm2_2_2}
S_0f=f\circ r,\text{ or }(S_0f)(x)=f(r(x)),\quad(x\in X).
\end{equation}
\item If $m_0\in L^\infty(X)$, we set
\[S_{m_0}=M(m_0)S_0,\]
or equivalently
\begin{equation}\label{eqfm2_2_3}
(S_{m_0}f)(x)=m_0(x)f(r(x)),\quad(x\in X,f\in L^\infty(X)).
\end{equation}
\item $r^{-1}(E):=\{\,x\in X\mid r(x)\in E\,\}$ for
$E\in\mathfrak{B}(X)$.
\[\mu\circ r^{-1}(E)=\mu(r^{-1}(E)),\quad(E\in\mathfrak{B}(X)).\]
\end{itemize}

\subsection{Properties of measures $\mu$ on $X$. Definitions}
\begin{enumerate}
\item {\it Invariance}:
\begin{equation}\label{eqfm4}
\mu\circ r^{-1}=\mu.
\end{equation}
\item {\it Strong invariance}:
\begin{equation}\label{eqfm5}
\int_Xf(x)\,d\mu=\int_X\frac{1}{\#r^{-1}(x)}\sum_{r(y)=x}f(y)\,d\mu\quad(f\in
L^\infty(X)).
\end{equation}
\item $W\colon X\rightarrow[0,\infty)$,
\begin{equation}\label{eqfm6}
(R_Wf)(x)=\sum_{r(y)=x}W(y)f(y).
\end{equation}
If $m_0\in L^\infty(X,\mu)$ is complex valued, we use the notation
$R_{m_0}:=R_W$ where $W(x)=|m_0(x)|^2/\#r^{-1}(r(x))$.
\begin{enumerate}
\item A function $h\colon X\rightarrow[0,\infty)$ is said to be an
eigenfunction for $R_W$ if
\begin{equation}\label{eqfm7}
R_Wh=h
\end{equation}
\item A Borel measure $\nu$ on $X$ is said to be a
left-eigenfunction for $R_W$ if
\begin{equation}\label{eqfm8}
\nu R_W=\nu,
\end{equation}
or equivalently
\[\int_XR_Wf\,d\nu=\int_Xf\,d\nu,\text{ for all }f\in L^\infty(X).\]
\end{enumerate}
\end{enumerate}

\begin{lemma}\label{lemfm1}
\begin{enumerate}
\raggedright
\item For measures $\mu$ on $X$ we have the implication
\textup{(\ref{eqfm5})}${}\Rightarrow{}$\textup{(\ref{eqfm4})}, but not conversely. 
\item
If $W$ is given and if $\nu$ and $h$ satisfy \textup{(\ref{eqfm8})} and
\textup{(\ref{eqfm7})} respectively, then
\begin{equation}\label{eqfm9}
d\mu:=h\,d\nu
\end{equation}
satisfies \textup{(\ref{eqfm4})}. 
\item If $\mu$ satisfies \textup{(\ref{eqfm5})},
and $m_0\in L^\infty(X)$, then $S_{m_0}$ is an isometry in
$L^2(X,h\,d\mu)$ if and only if
\[R_{m_0}h=h.\]
\end{enumerate}
\end{lemma}

\par
\subsection{Examples} We illustrate the definitions:
\begin{example}\label{exfm1}
Let $X=[0,1]=\br/\bz$. Fix $N\in\bz_+$, $N>1$. Let
\[r(x)=Nx\mod 1\]
\par
Invariance:
\begin{equation}\label{eqfme4}
\int_0^1f(Nx)\,d\mu(x)=\int_0^1f(x)\,d\mu(x)\quad(f\in
L^\infty(\br/\bz)).
\end{equation}
\par
Strong invariance:
\begin{equation}\label{eqfme5}
\frac{1}{N}\int_0^1\sum_{k=0}^{N-1}f\left(\frac{x+k}{N}\right)\,d\mu(x)=\int_0^1f(x)\,d\mu(x).
\end{equation}
\par
The Lebesgue measure $\mu=\lambda$ is the unique probability
measure on $[0,1]=\br/\bz$ which satisfies (\ref{eqfme5}).
\par
Examples of measures $\mu$ on $\br/\bz$ which satisfy
(\ref{eqfme4}) but not (\ref{eqfme5}) are
\begin{itemize}
\item $\mu=\delta_0$, the Dirac mass at $x=0$; \item
$\mu=\mu_{\bf{C}}$, the Cantor middle-third measure on $[0,1]$
(see \cite{DJ03}), i.e., $\mu_{\bf{C}}$ is determined by
\begin{itemize}
\item
$\frac{1}{2}\int\left(f\left(\frac{x}{3}\right)+f\left(\frac{x+2}{3}\right)\right)\,d\mu_{\bf{C}}(x)=\int
f(x)\,d\mu_{\bf{C}}(x),$ \item
$\mu_{\bf{C}}([0,1])=1,$ \item
$\mu_{\bf{C}}$ is supported on the middle-third Cantor set.
\end{itemize}
\end{itemize}
\end{example}

\begin{example}\label{exfme2}
Let $X=[0,1)=\br/\bz$, $\lambda$ the Lebesgue measure, $X_{\bf C}$
the middle-third Cantor set, $\mu_{\bf C}$ the Cantor measure.
\par
$r\colon X\rightarrow X$, $r(x)=3x\bmod 1$, $r_{\bf C}=r_{X_{\bf
C}}\colon X_{\bf C}\rightarrow X_{\bf C}$.
\par
Consider the following properties for a Borel probability measure
$\mu$ on $\br$:
\begin{equation}\label{eqfmeinv}
\int
f\,d\mu=\frac{1}{3}\int\left(f\Bigl(\,\frac{x}{3}\,\Bigr)+f\Bigl(\,\frac{x+1}{3}\,\Bigr)
+f\Bigl(\,\frac{x+2}{3}\,\Bigr)\right)\,d\mu(x);
\end{equation}
\begin{equation}\label{eqfmeinvc}
\int
f\,d\mu=\frac{1}{2}\int\left(f\Bigl(\,\frac{x}{3}\,\Bigr)
+f\Bigl(\,\frac{x+2}{3}\,\Bigr)\right)\,d\mu(x);
\end{equation}
Then (\ref{eqfmeinv}) has a unique solution $\mu=\lambda$.
Moreover (\ref{eqfmeinvc}) has a unique solution, $\mu=\mu_{\bf
C}$, and $\mu_{\bf C}$ is supported on the Cantor set $X_{\bf C}$.
\par
Let $\br/\bz=[0,1)$. Then $\#r^{-1}(x)=3$ for all $x\in[0,1)$. If
$x=x_1/3+x_2/3^2+\cdots$, $x_i\in\{0,1,2\}$, is the
representation of $x$ in base $3$, then $r(x)\sim(x_2,x_3,\dots)$,
and
$r^{-1}(x)=\{(0,x_1,x_2,\dots),(1,x_1,x_2,\dots),(2,x_1,x_2,\dots)\}$
\par
On the Cantor set $\#r_{\bf C}^{-1}(x)=2$ for all $x\in X_{\bf
C}$. If $x=x_1/3+x_2/3^2+\cdots$, $x_i\in\{0,2\}$, is
the usual representation of $X_{\bf C}$ in base $3$, then
$r_{\bf C}(x)=(x_2,x_3,\dots)$
and
$
X_{\bf C}\simeq\prod_{\mathbb{N}}\{0,2\}
$.
\par
In the representation $\prod_{\mathbb{N}}\bz_3$ of $X=[0,1)$,
$\mu=\lambda$ is the product (Bernoulli) measure with weights
$(1/3,1/3,1/3)$.
\par
In the representation $\prod_{\mathbb{N}}\{0,2\}$ of $X_{\bf C}$,
$\mu_{\bf C}$ is the product (Bernoulli) measure with weights
$(1/2,1/2)$.
\end{example}

\begin{example}\label{exfme3}
Let $N\in\bz_+$, $N\geq2$ and let $A=(a_{ij})_{i,j=1}^N$ be an $N$
by $N$ matrix with all $a_{ij}\in\{0,1\}$, and let
$X(A)$ be the corresponding state space
and let $r=r_A$.
\begin{lemma}\label{lemfme3_1}
Let $A$ be as above. Then
\[
\#r_A^{-1}(x)=\#\{\,y\in\{1,\dots,N\}\mid A(y,x_1)=1\,\}.
\]
\end{lemma}
\par
It follows that $r_A\colon X(A)\rightarrow X(A)$ is onto iff $A$ is {\it
irreducible}, i.e., iff for all $j\in\bz_N$, there exists an
$i\in\bz_N$ such that $A(i,j)=1$.
\end{example}
 The existence of strongly invariant measures is guaranteed by results of
 Ruelle \cite{Rue89} for the case of subshifts of finite type,
 Brolin \cite{Bro} for rational maps and
 \cite{Bal00} for a class of expanding maps.
\par
 One of the tools from operator theory which has been especially useful
 in the analysis of wavelets is multiplicity theory
for abelian $C^*$-algebras $\mathcal{A}$.

We first recall a few well known facts, see e.g., \cite{N}. By
Gelfand's theorem, every abelian $C^*$-algebra with unit is $C(X)$
for a compact Hausdorff space X; and every representation of
$\mathcal{A}$ is the orthogonal sum of cyclic representations.
While the cardinality of the set of cyclic components in this
decomposition is an invariant, the explicit determination of the
cyclic components is problematic, as the construction depends on
Zorn's lemma.
\par
Suppose now that $H$ is a Hilbert space with an isometry $S$ on it
and with a normal representation $\pi$ of $L^\infty(X)$ on $H$
that satisfies the covariance relation
\begin{equation}\label{eqmul2}
S\pi(g)=\pi(g\circ r)S,\quad(g\in L^\infty(X)).
\end{equation}
\par
Theorem \ref{th2_5} shows that there exists a Hilbert space $\hat
H$ containing $H$, a unitary $\hat S$ on $\hat H$  and a
representation $\hat\pi$ of $L^\infty(X)$ on $\hat H$ such that:
\[(V_n:=\hat S^{-n}(H))_n\text{ form an increasing sequence of
subspaces with dense union},\]
\[\hat S|_H=S,\]
\[\hat\pi|_H=\pi,\]
\[\hat S\hat\pi(g)=\hat\pi(g\circ r)\hat S.\]

\begin{theorem}\label{th2_5}
\textup{(i)} Let $H$ be a Hilbert space, $S$ an isometry on $H$. Then there
exist a Hilbert space $\hat H$ containing $H$ and a unitary $\hat
S$ on $\hat H$ such that
\begin{equation}\label{eq2_5_1}
\hat S|_{H}=S,
\end{equation}
\begin{equation}\label{eq2_5_2}
\cj{\bigcup_{n\geq0}\hat S^{-n}H}=\hat H.
\end{equation}
Moreover these are unique up to an intertwining isomorphism. \par
\textup{(ii)} If $\mathcal{A}$ is a $C^*$-algebra, $\alpha$ is an
endomorphism on $\mathcal{A}$ and $\pi$ is a representation of
$\mathcal{A}$ on $H$ such that
\begin{equation}\label{eq2_5_3}
S\pi(g)=\pi(\alpha(g))S\quad(g\in\mathcal{A});
\end{equation}
then there exists a unique representation $\hat\pi$ on $\hat H$
such that
\begin{equation}\label{eq2_5_4}
\hat\pi(g)|_H=\pi(g)\quad(g\in\mathcal{A}),
\end{equation}
\begin{equation}\label{eq2_5_5}
\hat S\hat\pi(g)=\hat\pi(\alpha(g))\hat S\quad(g\in\mathcal{A}).
\end{equation}
\end{theorem}

\begin{corollary}\label{cor2_6}
Let $X,r,$ and $\mu$ be as above. Let $I$ be a finite or countable
set. Suppose $H\colon X\rightarrow\mathcal{B}(l^2(I))$ has the property
that $H(x)\geq 0$ for almost every $x\in X$, and $H_{ij}\in
L^1(X)$ for all $i,j\in I$. Let
$M_0\colon X\rightarrow\mathcal{B}(l^2(I))$ such that
$x\mapsto\|M_0(x)\|$ is essentially bounded. Assume in addition
that
\begin{equation}\label{eq2_6_1}
\frac{1}{\#r^{-1}(x)}\sum_{r(y)=x}M_0^*(y)H(y)M_0(y)=H(x),\text{ for a.e. }x\in X.
\end{equation}
Then there exists a Hilbert space $\hat K$, a unitary operator
$\hat U$ on $\hat K$, a representation $\hat\pi$ of $L^\infty(X)$
on $\hat K$, and a family of vectors $(\varphi_i)\in\hat K$, such
that:
\[\hat U\hat\pi(g)\hat U^{-1}=\hat\pi(g\circ r)\quad(g\in L^\infty(X)),\]
\[\hat U\varphi_i=\sum_{j\in I}\hat\pi((M_0)_{ji})\varphi_j\quad(i\in I),\]
\[\ip{\varphi_i}{\hat\pi(f)\varphi_j}=\int_XfH_{ij}\,d\mu\quad(i,j\in I,\;f\in L^\infty(X)),\]
\[
\cj{\operatorname*{span}}\{\,\hat\pi(f)\varphi_i\mid n\geq0,\;f\in L^\infty(X),\;i\in I\,\}=\hat K.
\]
These are unique up to an intertwining unitary isomorphism.
\textup{(}All functions are assumed weakly measurable in the sense that 
$x\mapsto\ip{\xi}{F(x)\eta}$ is measurable for all $\xi,\eta\in l^2(I)$.\textup{)}
\end{corollary}

Note that $\hat S$ maps $V_1$ to $V_0$, and the covariance
relation implies that the representation $\hat\pi$ on $V_1$ is
isomorphic to the representation $\pi^r\colon g\mapsto\pi(g\circ r)$ on
$V_0$. Therefore we have to compute the multiplicity of the
latter, which we denote by $d^r_{V_0}$.
\par
By the spectral theorem there exists a unitary isomorphism
$J\colon H(=V_0)\rightarrow L^2(X,d_{V_0},\mu)$, where, for a {\it
multiplicity function} $d\colon X\rightarrow\{0,1,\dots,\infty\}$, we use
the notation:
\[
L^2(X,d,\mu):=\left\{\,f\colon X\rightarrow\cup_{x\in X}\bc^{d(x)}\Bigm|
f(x)\in\bc^{d(x)},\;\int_X\|f(x)\|^2\,d\mu(x)<\infty\,\right\}.
\]
In
addition $J$ intertwines $\pi$ with the representation of
$L^\infty(X)$ by multiplication operators, i.e.,
\[(J\pi(g)J^{-1}(f))(x)=g(x)f(x)\,\quad(g\in L^\infty(X),\;f\in
L^2(X,d_{V_0},\mu),\;x\in X).\]

\begin{remark}
Here we are identifying $H$ with $L^2(X,d_{V_0},\mu)$ via the {\it
spectral representation}. We recall the details of this
representation $H\ni f\mapsto\tilde f\in L^2(X,d_{V_0},\mu)$.
\par
Recall that any normal representation $\pi\in Rep(L^\infty(X),H)$
is the orthogonal sum
\begin{equation}\label{eqmulrem*}
H=
\sideset{}{^{\smash{\oplus}}}{\sum}\limits_{k\in C}
\left[\pi(L^\infty(X))k\right],
\end{equation}
where the set $C$ of vectors $k\in H$ is chosen such that
\begin{itemize}
\item $\quad\quad\quad\quad\quad\quad\quad\quad\quad\quad\|k\|=1,$
\begin{equation}\label{eqmulrem**}
\ip{k}{\pi(g)k}=\int_Xg(x)v_k(x)^2\,d\mu(x),\text{ for all }k\in
C;
\end{equation}
\item $\ip{k'}{\pi(g)k}=0,\quad g\in L^\infty(X),\;k,k'\in C,k\neq
k';\text{ orthogonality}.$ \end{itemize}
\par
The formula (\ref{eqmulrem*}) is obtained by a use of Zorn's
lemma. Here, $v_k^2$ is the Radon-Nikodym derivative of
$\ip{k}{\pi(\,\cdot\,)k}$ with respect to $\mu$, and we use that $\pi$
is assumed normal.
\par
For $f\in H$, set
\[f=
\sideset{}{^{\smash{\oplus}}}{\sum}\limits_{k\in C}
\pi(g_k)k,\quad g_k\in L^\infty(X)\]
and
\[\tilde f=
\sideset{}{^{\smash{\oplus}}}{\sum}\limits_{k\in C}
g_kv_k\in L_\mu^2(X,l^2(C)).\]
Then $Wf=\tilde f$ is the desired spectral transform, i.e.,
\[W\text{ is unitary},\]
\[W\pi(g)=M(g)W,\]
and
\[\|\tilde f(x)\|^2=\sum_{k\in C}|g_k(x)v_k(x)|^2.\]
Indeed, we have
\[\int_X\|\tilde f(x)\|^2\,d\mu(x)=\int_X\sum_{k\in
C}|g_k(x)|^2v_k(x)^2\,d\mu(x)=\sum_{k\in
C}\int_X|g_k|^2v_k^2\,d\mu.\]
\[=\sum_{k\in C}\ip{k}{\pi(|g_k|^2)k}=\sum_{k\in
C}\|\pi(g_k)k\|^2=\left\|
\sideset{}{^{\smash{\oplus}}}{\sum}\limits_{k\in C}
\pi(g_k)k\right\|^2_H=\|f\|^2_{H}.\]
\par
It follows in particular that the multiplicity function
$d(x)=d_{H}(x)$ is
\[
d(x)=\#\{\,k\in C\mid v_k(x)\neq0\,\}.
\]
Setting
\[
X_i:=\{\,x\in X\mid d(x)\geq i\,\}\quad (i\geq 1),
\]
we see that
\[H\simeq
\sideset{}{^{\smash{\oplus}}}{\sum}
L^2(X_i,\mu)\simeq L^2(X,d,\mu),\]
and the isomorphism intertwines $\pi(g)$ with multiplication
operators.
\end{remark}

\begin{theorem}\label{thmul1}
\textup{(i)} $V_1=\hat S^{-1}(H)$ is invariant for the representation
$\hat\pi$. The multiplicity functions of the representation
$\hat\pi$ on $V_1$, and on $V_0=H$, are related by
\begin{equation}\label{eqmul3}
d_{V_1}(x)=\sum_{r(y)=x}d_{V_0}(y)\quad(x\in X).
\end{equation}
\par
\textup{(ii)} If $W_0:=V_1\ominus V_0=\hat S^{-1}H\ominus H$, then
\begin{equation}\label{eqmul4}
d_{V_0}(x)+d_{W_0}(x)=\sum_{r(y)=x}d_{V_0}(y)\quad(x\in X).
\end{equation}
\end{theorem}

With the spaces $H_n$ in (\ref{eqhn}), our global Hilbert
space of discrete $L^2$-martingales thus acquires a concrete form, and we are
now able to present our multiresolution/wavelet result below.

\begin{theorem}\label{thpr_1}
If $h\in L^1(X)$, $h\geq0$ and $Rh=h$, then there exists a unique
measure $\hat\mu$ on $X_\infty$ such that
\[\hat\mu\circ\theta_n^{-1}=\omega_n\quad(n\geq0),\]
where \begin{equation}\label{eqpr_6}
\omega_n(f)=\int_XR^n(fh)\,d\mu\quad(f\in L^\infty(X)).
\end{equation}
\end{theorem}
\par
We give now a different representation of
the construction of the covariant system associated to $m_0$ and
$h$.
\par
Let
\begin{equation}\label{eqhn}
H_n:=\{\,f\in L^2(X_\infty,\hat\mu)\mid
f=\xi\circ\theta_n,\;\xi\in L^2(X,\omega_n)\,\}.
\end{equation}
Then $H_n$ form an
increasing sequence of closed subspaces which have dense union.
\par
We can identify the functions in $H_n$ with functions in
$L^2(X,\omega_n)$, by
\[i_n(\xi)=\xi\circ\theta_n\quad(\xi\in L^2(X,\omega_n)).\]
The definition of $\hat\mu$ makes $i_n$ an isomorphism between
$H_n$ and $L^2(X,\omega_n).$
\par
Define
\[
\mathcal{H}:=\left\{(\xi_0,\xi_1,\dots)\Bigm| \xi_n\in
L^2(X,\omega_n),\,R(\xi_{n+1}h)=\xi_nh,\,
\sup_n\!\int_X\!R^n(|\xi_n|^2h)\,d\mu<\infty\right\},
\]
with the scalar
product
\[\ip{(\xi_0,\xi_1,\dots)}{(\eta_0,\eta_1,\dots)}=\lim_{n\rightarrow\infty}\int_XR^n(\cj\xi_n\eta_nh)\,d\mu.\]

\begin{theorem}\label{thpr_5}
The map $\Phi\colon L^2(X_\infty,\hat\mu)\rightarrow\mathcal{H}$ defined
by
\[\Phi(f)=(i_n^{-1}(P_nf))_{n\geq0},\]
where $P_n$ is the projection onto $H_n$ is an isomorphism.
Moreover this transform $\Phi$ satisfies the following three properties,
and it is determined by them:
\[\Phi U\Phi^{-1}(\xi_n)_{n\geq0}=(m_0\circ
r^n\,\xi_{n+1})_{n\geq0},\]
\[\Phi\pi(g)\Phi^{-1}(\xi_n)_{n\geq0}=(g\circ
r^n\,\xi_n)_{n\geq0},\]
\[\Phi\varphi=(1,1,\dots).\]
\end{theorem}
\begin{acknowledgements}
The authors thank Dick Gundy for helpful discussions,
and the National Science Foundation for partial support.
\end{acknowledgements}

\end{document}